\documentclass{article}

\usepackage{amsmath, amssymb, eucal}

\begin{document}

\righthyphenmin=2
\hfuzz=0.5pt

\def\R{{\Bbb R}}
\def\C{{\Bbb C}}
\def\H{{\Bbb H}}
\def\K{{\Bbb K}}

\def\le{{\leqslant}}
\def\ge{{\geqslant}}

\def\cR{{\cal R}}
\def\fr{{\frak r}}
\def\fP{{\frak P}}
\def\cL{{\cal L}}
\def\kappa{{\varkappa}}

\def\GL{{\rm GL}}
\def\U{{\rm U}}
\def\OO{{\rm O}}
\def\Sp{{\rm Sp}}

\def\ov{\overline}

\def\Herm{{\rm Herm}}
\def\Mat{{\rm Mat}}
\def\Pos{{\rm Pos}}

\def\res{\mathop{\rm res}}

\def\Re{\mathop{\rm Re}}

\def\tr{\mathop{\rm tr}}

\newcounter{sec}
\renewcommand{\theequation}{\arabic{sec}.\arabic{equation}}

\begin{center}

{\Large\bf
Rayleigh triangles and non-matrix interpolation
of matrix beta-integrals
}

\bigskip

{\sc Yu. A. Neretin}\footnote{
Supported in part by  NWO Grant 047-008-009}

\end{center}

{\small\begin{quotation}
A Rayleigh triangle of size $n$
is a collection of
$n(n+1)/2$  real numbers $\lambda_{kl}$,
where $1\le l\le k\le n$; these numbers decrease
with growth of
$k$ for fixed  $l$ and increase with growth of
$k$ for fixed $k-l$. We construct a family of
beta-integrals over the space of Rayleigh triangles;
these integrals interpolate the matrix beta-integrals of
the Siegel,
Hua Loo Keng and Gindikin types with respect to the dimension
of the basic field
 ($\R$, $\C$ or the quaternions $\H$).
We also interpolate the Hua--Pickrell measures
on the inverse limits of the symmetric spaces
  $\U(n)$, $\U(n)/\OO(n)$,
$\U(2n)/\Sp(n)$.

Our family of integrals also includes
the Selberg integral.
\end{quotation} }

{\bf 0.1. Interpolation of homogeneous spaces.}
It is well known that the representation theory
of semisimple groups admits a partial interpolation
on the level of special functions.
The most old (and the most simple) construction
of this kind is the Hankel transform
  \cite{Han}, \cite{Wat}.
Let $f$ be a function on  $\R^n$
depending only on  $r=(x_1^2+\dots+x_n^2)^{1/2}$.
Consider its Fourier transform
$$
\widehat f(a)=
\frac 1{(2\pi)^{n/2}}\int_{\R^n}
f(r)e^{-iax} dx
.
$$
Averaging $e^{iax}$ over spheres of radius
$\rho=(a_1^2+\dots+a_n^2)^{1/2}$,
we obtain that
 $\widehat f(\rho)$
is given by
\begin{equation}
\widehat f(\rho)=
H_n f(\rho):=
\int_0^\infty
J_{n/2-1}(r\rho)
f(r) r^{n/2} dr
,
\end{equation}
where
$J_\nu(z)$ is a Bessel function
$$J_\nu(z)=
\sum_{m=0}^\infty \frac{ (-1)^m \left(\frac z2 \right)^{2m+\nu}}
                   {m!\,\Gamma(m+\nu+1)}
.
$$

Obviously, the operator
 $H_n$ is a unitary operator
in the space
$L^2(\R_+, r^{n-1}dr)$
and satisfies to the property
$H_n^2=1$.

By our construction,
the number $n$ is integer,
since $n$ is the dimension  of the space
 $\R^n$. Nevertheless the  {\it Hankel transform}%
\footnote{We use a nonstandard notation for the Hankel transform.}
(0.1) is well defined for all complex
 $n$ satisfying the condition $\Re n>0$.
By the Hankel theorem,
for all $n$ the identity  $H_n^2=1$ holds; moreover
for real  $n>0$ the operator $H_n$ is unitary
in the space $L^2(\R_+, r^{n-1}dr)$.

Similarly, the index hypergeometric transform
for some exceptional values of its parameters gives
the spherical transforms for all the hyperbolic
symmetric spaces
(\cite{Ole}, \cite{Koo}, \cite{Neri}).
For semisimple groups of rang $>1$
the interpolation can be performed only for functions
depending on eigenvalues or singular values of matrices.
Various topics related to this interpolation are discussed
\cite{HO}, \cite{Hec},  \cite{Mac}, \cite{Che}, \cite{CO}.

\smallskip

{\bf 0.2. Purposes of the paper.}
This paper has two aims:

--- to construct an interpolation of the Hua--Pickrell
measures (\cite{Pic},
\cite{Nerh}, \cite{BO}, \cite{Ols})
on the inverse limits of symmetric spaces.

--- to construct an interpolation for matrix
  $B$-functions (see \cite{Gin}, \cite{Nerb});
in this case we intend to replace the dimension of
the ground field
$\dim \K$, where $\K=\R$, $\C$ or quaternionic algebra $\H$,
by an arbitrary complex number.

 Recall that the second problem is not standard
from the point of view of the modern
interpolation ideology, since the
matrix
  $B$-function is not a function of matrix eigenvalues.
The both goals are achieved only partially.
Namely this is done for  simple groups of the  series
 $A$.

Emphasis that our family of beta-integrals
includes the Selberg integral.

\smallskip

{\bf 0.3.  Hua--Pickrell measures and their interpolation.}
Now we formulate in more details the problem
about the Hua--Pickrell measures.
Below we do not return to this topic and the rest of our text
is a 'play in integrals'.

By $\U(n)$ we denote the group of unitary matrices
of order $n$.
Consider the map $\Upsilon: \U(n+1)\to\U(n)$,
given by
$$\Upsilon
\begin{pmatrix} P&Q\\R&T\end{pmatrix}
=T-R(1+P)^{-1}Q
,$$
where
$\begin{pmatrix} P&Q\\R&T\end{pmatrix}\in\U(n+1)$
is a block matrix of size  $(1+n)\times (1+n)$.
It can be easily checked that
the image of the Haar measure on $\U(n+1)$
 under the map
$\Upsilon$ is the Haar measure on $\U(n)$.
Hence we can construct the inverse limit
$\frak U$ of the chain
\begin{equation}
\dots
\stackrel{\Upsilon}{\longleftarrow}
\U(n-1)
\stackrel{\Upsilon}{\longleftarrow}
\U(n)
\stackrel{\Upsilon}{\longleftarrow}
\U(n+1)
\stackrel{\Upsilon}{\longleftarrow}
\dots
.
\end{equation}
Consider a sequence of matrices
 $g_n\in\U(n)$ satisfying the condition
$$
g_n=\Upsilon (g_{n+1}    )
$$
(the set of such chains is in an one-to-one correspondence
with the space $\frak U$).
Denote by $\nu_{n1},\dots,\nu_{nn}$
the eigenvalues of the matrix
$g_n$.

Thus, for each element of the space
 $\frak U$ and each $n=1,2,\dots$
we assign
the collection of points
$\nu_{n1},\dots,\nu_{nn}$   of the circle $|\nu|=1$.

In other words, we obtain a stochastic process with a discrete
time, its value in a moment $n$ is an $n$-point
subset of the circle.

The construction admits three following variations,
see \cite{Nerh}.

First, we can replace the Haar measure
$dg$ на $\U(n)$ by the two-parameter family
of the Hua--Pickrell measures
$$
c_n(\sigma)\det(1+g)^{\sigma} \det(1+\overline g)^{\overline\sigma} dg
,$$
where $\sigma\in\C$ is fixed
and the normalizing factor $c_n(\sigma)$
is defined by the condition: the measure of the whole group
 $\U(n)$ is 1.
These measures form a consistent system with respect
to the maps
$\Upsilon$,
and we again can consider the projective
limit of our chain.

Second, the Hua--Pickrell measures can be
replaced by a family of measures depending
on an infinite family of complex parameters
$\sigma_1,\sigma_2, \dots$. In this case,
the measure on the group $\U(n)$ is given by
\begin{multline*}
c_n(\sigma_1,\sigma_2,\dots)
\det(1+g)^{\sigma_n} \det(1+\overline g)^{\overline\sigma_n}
\times\\ \times
\prod\limits_{j=1}^{n-1} \det(1+[g]_j)^{\sigma_j-\sigma_{j+1}}
\prod\limits_{j=1}^{n-1}
\det(1+[\overline g]_j)^{\overline\sigma_j-\overline\sigma_{j+1}}
dg
,
\end{multline*}
where $[g]_j$ denotes the left upper $j\times j$
block of the matrix   $g$.

Thirdly, this construction survives for all the classical
compact symmetric spaces.
Consider the inverse limits
${\frak U}/{\frak O}$ and
${\frak U}/{\frak S \frak p}$
of two chains
\begin{gather}
\dots\longleftarrow \U(n-1)/\OO(n-1)
             \longleftarrow \U(n)/\OO(n)
             \longleftarrow         \dots           \\
\dots\longleftarrow \U(2n-2)/\Sp(2n-2)
             \longleftarrow \U(2n)/\Sp(2n)
             \longleftarrow  \dots
\end{gather}

Apparently, the three spaces
$\frak U$, ${\frak U}/{\frak O}$, ${\frak U}/{\frak S\frak p}$
themselves can not be interpolated.
However our basic construction (Theorem 3.1)
can be interpreted in the following way.
We construct a large class of Markov stochastic processes,
whose value in a moment $n$ is an $n$-point subset
on the circle.
In particular, this class contains
processes on the circle related to the chains
 (0.2) -- (0.4).

\smallskip

I am grateful to G.I.Olshanski for discussion of this subject.

    \medskip

{\large\bf \S\,1. Main results}

\addtocounter{sec}{1}
\setcounter{equation}{0}

\medskip

{\bf 1.1. Rayleigh triangles.}
Let $A$ be an Hermitian $n\times n$ matrix,
let $B$ be its left upper
 $(n-1)\times (n-1)$ block.
Denote by
$\mu_1\le \mu_2\le\dots \le \mu_n$
the eigenvalues of
the matrix $A$
and by
 $\lambda_1\le\lambda_2\le\dots \le\lambda_{n-1}$
the eigenvalues of the matrix $B$.

By the Rayleigh theorem (see for instance \cite{Bha}, \cite{Gan}),
\begin{equation}
\mu_1\le\lambda_1\le \mu_2\le\lambda_2\le\mu_3
\le \dots \le  \lambda_{n-1}\le\mu_n
.
\end{equation}

We define  a {\it  Rayleigh triangle} $\cL$ of size $n$
as a collection of
 $n(n+1)/2$ real numbers
\begin{equation}
\begin{array} {ccccccccccc}
 & & & &   & \lambda_{11}  &  & & & &  \\
 & & & &  \lambda_{21}  & & \lambda_{22}  & & & &  \\
 & & &\lambda_{31} &   &\lambda_{32} &  &\lambda_{33} & & &  \\
 & &\cdot &\cdot & \cdot  &\cdot & \cdot &\cdot &\cdot & &  \\
 &\cdot &\cdot &\cdot & \cdot  &\cdot & \cdot &\cdot &\cdot &\cdot &  \\
\lambda_{n1} & &\lambda_{n2} & &\lambda_{n3}
      &\dots & \lambda_{n(n-2)} & &\lambda_{n(n-1)}
        &  & \lambda_{nn}
,
\end{array}
\end{equation}
satisfying the system of inequalities
\begin{equation}
\dots\le \lambda_{(j+1)k}\le \lambda_{jk}\le \lambda_{(j+1)(k+1)}\le \dots
\end{equation}
for all the possible $j$, $k$.
In other words, the numbers $\lambda_{j\alpha}$
decrease  south--west--ward direction
and increase  south--east--ward.

We denote the space
 of all the Rayleigh triangles of the size
$n$ by $\cR_n$.
By $\cR_n[a,b]$ we denote the set of all the
Rayleigh triangles of size
$n$ satisfying the condition
$$a\le \lambda^{(n)}_1 , \qquad \lambda^{(n)}_n \le b$$
(hence all the numbers $\lambda_{j\alpha}$
lie between $a$ and $b$).

It is known, that some matrix integrals can be reduced
to the integration over the spaces
$\cR_n[a,b]$;
at least  this was a way of the evaluation
of spherical functions of the group $\GL(n,\C)$
in Gelfand's and Naimark's book \cite{GN}.
A.Okounkov and G.I.Olshanski \cite{OO}
obtained a representation of the Jack polynomials
as some integrals over the space of Rayleigh triangles.
This work was a standpoint for the  present paper;
 see also  Kazarnovsky-Krol's work
\cite{Kaz1}.

\smallskip

{\bf 1.2. Natural measures on $\cR_n$.}
By $\K$ we denote $\R$, $\C$ or the quaternionic algebra
$\H$.
Assume
$$
\theta=\dim\K/2
.
$$
Let us define a natural measure  $\{dx\}$
on  $\K$:
for the real numbers it is the Lebesgue measure,
for the complex numbers $x=u+iv$ we put
 $\{d(u+iv)\}=du\,dv$,
and for the quaternions  $x=u+iv+jw+kz$
 we give $\{d(u+iv+jw+kz\}=du\,dv\,dw\,dz$.

\smallskip

Consider the space
$\Herm_\K(n)$
of all the Hermitian matrices%
\footnote{A square matrix $Z$ over $\K$ is Hermitian
if its matrix elements $z_{\alpha\beta}$
satisfy the condition $z_{\alpha\beta}=\overline z_{\beta\alpha}$.}
  of size  $n\times n$
 over $\K$.
Equip this space by the Lebesgue measure
$$
\prod_{k>l} \{dt_{kl}\} \prod_k dt_{kk}
.$$

For a matrix $T$ denote by $[T]_j$
its left upper block of size
 $j\times j$.

For each block $[T]_j$ of the matrix
$T$ denote by
$$\lambda_{j1}\le \lambda_{j2} \le\dots\le \lambda_{jj}$$
the eigenvalues of the matrix $[T]_j$.
Thus we obtain the map
$$
\Herm_\K(n)\to \cR_n(-\infty,\infty)
.$$

\smallskip

{\sc Proposition 1.1.}\footnote%
{The author does not know is this statement published
somewhere; it can be observed from the calculation
\cite{GN}, \S II.9. It is known in a folklore
(since the integral representation of the Jack multivariate
hypergeometric functions \cite{OO}, \cite{Kaz1}
ideologically depend on this  statement.}
{\it  The image of the Lebesgue measure on
$\Herm_\K(n)$ under the map
$\Herm_\K(n)\to \cR_n$ is
\begin{equation}
C_n(\theta)\fr_\theta(\cL)\,d \cL
,
\end{equation}
 where
\begin{multline}
\fr_\theta(\cL) = \frac
{\prod\limits_{2\le j\le n} \quad
\prod\limits_{1\le \alpha\le j-1,\,\,\,
                               1\le p\le j}
|\lambda_{(j-1)\alpha}-\lambda_{jp}|^{\theta-1} }
{\prod\limits_{2\le j\le n-1}\quad
\prod\limits_{1\le \alpha<\beta\le j}
  (\lambda_{j\beta}-\lambda_{j\alpha})^{2\theta-2}}
\times \\ \times
\prod\limits_{1\le p<q\le n}
(\lambda_{nq} -\lambda_{np})
,
\end{multline}
here $d \cL$ is the Lebesgue measure
on the space of  Rayleigh triangles
$$
d\cL=
 \prod\limits_{1\le j\le n}\quad \prod_{1\le\alpha\le j}
  d\lambda_{j\alpha}
$$
and $C_n(\theta)$ is the normalizing constant
}
$$
C_n(\theta)=\frac{\pi^{n(n-1)\theta}}
{\Gamma^{n(n-1)/2} (\theta)}
.
$$

\smallskip

By $\Mat_\K(n,m)$ we denote the space
of matrices of  size $n\times m$ over
$\K$. Assume $n\le m$,
i.e., the number of rows  does not exceed the number
of columns.
By $\{T\}_j$ we denote the matrix consisting
of the first
 $j$ rows of the matrix $T$.
 For each
$j=1,2,\dots,n$, we denote by
$$\lambda_{j1}\le \lambda_{j2} \le\dots\le \lambda_{jj}$$
the eigenvalues of the matrix
$\{T\}_j\{T\}_j^*$,
obviously all these numbers
 are nonnegative.
Thus we obtain the map
\begin{equation}
\Mat_\K(n,m)\to\cR_n[0,\infty)
.
\end{equation}

\smallskip

{\sc Proposition 1.2.}
{\it   The image of the Lebesgue measure
on $\Mat(n,m)$ under the map
{\rm (1.6)} is
$$
\widetilde C_{nm}(\theta)\psi(\cL)\fr_\theta(\cL)\,d \cL
,
$$
where $\fr_\theta(\cL)$ is given by {\rm (1.5)},
$$
\psi(\cL)=
\prod\limits_{p=1}^n \lambda_{n p}^{(m-n+1)\theta-1}
,
$$
and             }
$$
\widetilde C_{nm}(\theta)=\frac{\pi^{mn\theta}}
  {\Gamma^{n(n-1)}(\theta)\prod_{j=1}^n\Gamma(
(m-j+1)\theta)}
.
$$

{\bf 1.3. Ultra-beta-integrals.}
We deduce two following integrals
(1.7)--(1.8)
over the space of Rayleigh triangles.

Fix complex numbers  $\sigma_j$, $\tau_j$,
where $j=1,\dots, n$,
and also the complex numbers $\theta_{j\alpha}$,
where $j$, $\alpha$ range in the domain
$1\le\alpha\le j\le n-1$. Then
\begin{multline}
\int\limits_{\cR_n[0,\infty)}
\prod\limits_{j=1}^{n-1}\prod\limits_{\alpha=1}^j
\frac{\lambda_{j\alpha}^{\sigma_j-\sigma_{j+1}-\theta_{j\alpha}} }
 {(1+\lambda_{j\alpha})^{\tau_j-\tau_{j+1}+\theta_{j\alpha}} }
\cdot \prod\limits_{p=1}^n
  \frac {\lambda_{np}^{\sigma_n-1}}
     {(1+\lambda_{np})^{\tau_n}}
\times\\ \times
\prod\limits_{j=1}^{n-1}
\frac{
  \prod\limits_{1\le\alpha\le j,\, 1\le p\le j+1}
   |\lambda_{j\alpha}-\lambda_{(j+1)p}|^{\theta_{j\alpha}-1}
   }
   {
  \prod\limits_{1\le\alpha<\beta \le j}
  (\lambda_{j\beta}-\lambda_{j\alpha})^{\theta_{j\alpha}+\theta_{j\beta}-2}
  }
\prod\limits_{1\le p< q\le n} (\lambda_{nq}-\lambda_{np})\, d\cL=
\\=
\prod\limits_{1\le\alpha\le j\le n-1} \Gamma(\theta_{j\alpha})
\cdot
\prod\limits_{j=1}^n
\frac
  {\Gamma(\sigma_j)
    \Gamma(\tau_j-\sigma_j-\sum_{\alpha=1}^{j-1} \theta_{(j-1)\alpha})}
{\Gamma(\tau_j)}
.
\end{multline}
The conditions of absolute convergence
of this integral are
$$\Re\theta_{j\alpha}>0,\quad \Re\sigma_j>0,\quad
\Re\tau_j>\Re\sigma_j+\sum_{\alpha=1}^{j-1} \Re\theta_{(j-1)\alpha}
.
$$

As above, let
 $\tau_j$, $\sigma_j$, $\theta_{j\alpha}$
be complex numbers. Then
\begin{multline}
\int\limits_{\cR_n(-\infty,\infty)}
\prod\limits_{j=1}^{n-1}\prod\limits_{\alpha=1}^j
(1+i\lambda_{j\alpha})^{-\sigma_j+\sigma_{j+1}-\theta_{j\alpha}}
(1-i\lambda_{j\alpha})^{-\tau_j+\tau_{j+1}-\theta_{j\alpha}}
\times\\ \times \prod\limits_{p=1}^n
(1+i\lambda_{np})^{-\sigma_n}
(1-i\lambda_{np})^{-\tau_n}
\times\\ \times
\prod\limits_{j=1}^{n-1}
\frac{
  \prod\limits_{1\le\alpha\le j,\, 1\le p\le j+1}
   |\lambda_{j\alpha}-\lambda_{(j+1)p}|^{\theta_{j\alpha}-1}
   }
   {
  \prod\limits_{1\le\alpha<\beta \le j}
  (\lambda_{j\beta}-\lambda_{j\alpha})^{\theta_{j\alpha}+\theta_{j\beta}-2}
  }
\prod\limits_{1\le p< q\le n} (\lambda_{nq}-\lambda_{np})\, d\cL
=\\=
\pi^{n} 2^{2n-\sum\limits_{j=1}^n(\sigma_j+\tau_j)}
\prod\limits_{1\le\alpha\le j\le n-1} \Gamma(\theta_{j\alpha})
\cdot
\prod\limits_{j=1}^n
\frac{\Gamma(\sigma_j+\tau_j-1-\sum_{\alpha=1}^{j-1} \theta_{(j-1)\alpha})}
    {\Gamma(\sigma_j)\Gamma(\tau_j)}
,
\end{multline}
where $i$ is the imaginary unit.
The conditions of the absolute convergence of this integral
is
$$
\Re\theta_{j\alpha}>0; \quad
\Re\bigl\{\sigma_j+\tau_j-1-
\sum_{\alpha=1}^{j-1} \theta_{(j-1)\alpha})\bigr\}
>0
.
$$
In \S 3 we derive three more integrals
(3.1)--(3.3)  in the same spirit.

\smallskip

{\bf 1.4. Matrix integrals.}
Ultra-beta-integrals include many various
matrix beta- and gamma-integrals
of Whishart--Ingham
\cite{Ing},
Siegel \cite{Sie}, Hua Loo Keng \cite{Hua}
and Gindikin \cite{Gin} types.%
\footnote{The integrals from  author's work
\cite{Nerb} are not absorbed by ultra-beta-integrals.}
We present two series of such integrals now
and some more series in п.4.3.
I have not seen some of these integrals
in the literature.
Apparently, varying symmetry conditions
on matrices it is possible to obtain
many other integrals in the same style.

Consider the space $\Pos_\K(n)$
of Hermitian positive definite matrices $T$
of size  $n\times n$ over field $\K=\R,\C,\H$;
as above, let $\theta=\frac 12\dim\K$.
Fix complex numbers
$\sigma_1,\dots, \sigma_n$
and $\tau_1,\dots,\tau_n$.
The beta-integrals of Gindikin
(see \cite{Gin}, see also {FK},
they extend the Siegel integrals \cite{Sie})
have the form
\begin{multline*}
\int\limits_{\Pos_\K(n)}
\prod\limits_{j=1}^{n-1}
\frac{\det^{\sigma_j-\sigma_{j+1}-\theta}[T]_j}
     {\det^{\tau_j-\tau_{j+1}+\theta}[1+T]_j}
\frac{
\det^{\sigma_n-1} T }
{\det^{\tau_n}(1+T)  }
dT
 =  \\ =
         \pi^{n(n-1)/2}
\prod\limits_{j=1}^n
 \frac{\Gamma(\sigma_j)\Gamma(\tau_j-\sigma_j-(j-1)\theta)}
       {\Gamma(\tau_j)}
.
\end{multline*}

Let us explain how to reduce them to our integrals.
In notation of Subsection 1.4,
the measure $dT$ can be rewritten in the form (1.4),
the integrand transforms to
$$
\prod\limits_{j=1}^n \frac{\lambda_{j\alpha}^{\sigma_j-\sigma_{j+1}-\theta}}
                   {(1+\lambda_{j\alpha})^{\tau_j-\tau_{j+1}+\theta}}
\prod\limits_{p=1}^n \frac{\lambda_{np}^{\sigma_n-1}}
                   {(1+\lambda_{np})^{\tau_n}}
,$$
and we obtain a partial case of the integral
(1.7).

Similarly,
\begin{multline*}
\int\limits_{\Herm_n(\K)}
\prod\limits_{j=1}^{n-1}
\det(1+i[T]_j)^{-\sigma_j+\sigma_{j+1}-\theta}
\det(1-i[T]_j)^{-\tau_j+\tau_{j+1}-\theta}
\times \\ \times
 \det (1+iT)^{-\sigma_n}\det(1-iT)^{-\tau_n}
dT
=\\ =
\pi^{n(n+1)/2} 2^{2n-\sum \sigma_j -\sum \tau_j}
\prod_{j=1}^n
\frac{\Gamma(\tau_j+\sigma_j-1-(j-1)\theta)}
  {\Gamma(\sigma_j)\Gamma(\tau_j)}
.
\end{multline*}

For the field
 $\K=\C$ this integral is evaluated in
\cite{Ners}, \cite{Nerh},
one of partial cases of this integral
is the Hua Loo Keng integral \cite{Hua}
$$
\int\limits_{\Herm_n(\K)}
\det(1+ T^2)^{-\sigma} dT
,$$
it corresponds to
$$
\sigma_j=\tau_j=\sigma-(n-j)\theta
.
$$

Further,
\begin{multline*}
\int\limits_{\Mat_{m,n}(\K)}
\prod\limits_{j=1}^{n-1}
\frac {\det(\{T\}_j \{T\}_j^* )^{\sigma_j-\sigma_{j+1}-\theta}}
{\det(1+\{T\}_j \{T\}_j^* )^{\tau_j-\tau_{j+1}+\theta}}
\cdot \frac{\det(T^*T)^{\sigma_n-(m-n+1)\theta}} {\det(1+TT^*)^{\tau_n}}
dT
= \\ =
\pi^{mn\theta}
\prod\limits_{j=1}^n \frac{\Gamma(\sigma_j)
   \Gamma(\tau_j-\sigma_j-(j-1)\theta)}
{\Gamma(\tau_j)\Gamma((m-j+1)\theta)}
.
\end{multline*}
This series of integrals includes
the Hua Loo Keng integrals
\cite{Hua}
$$
\int\limits_{\Mat_{m,n}(\K)}
\det(1+ T^2)^{-\tau} dT
.
$$
Some more series of matrix integrals
are discussed below in
4.3.

\smallskip

{\bf 1.5. Projective systems of measures.}
Fix an infinite sequence of complex numbers
$\sigma_1,\sigma_2,\dots$. Assume  $\tau_j=\ov\sigma_j$.
Fix also reals numbers $\theta_{j\alpha}$,
where $j$, $\alpha$ range in $1\le\alpha\le j$.
Let these numbers satisfy
the inequalities
$$
\theta_{j\alpha}>0; \qquad
2\Re\sigma_j>1+\sum_{\alpha=1}^{j-1} \theta_{(j-1)\alpha}
.
$$

Further, fix $n$. Denote by
$\fP_n(\cL)$ the integrand in (1.8).
Consider the probability measure
$$\kappa_n=c_n\cdot \fP_n(\cL) d\cL$$
on  $\cR_n(-\infty, \infty)$,
where $c_n$ is a normalizing constant
(we choice it to make the measure probabilitic).

Consider the map
$$
\Upsilon:\cR_{n+1}(-\infty,\infty)\to
\cR_{n}(-\infty,\infty)
,
$$
that erases the last row of a Rayleigh triangle
$\cL\in\cR_{n+1}(-\infty,\infty)$.

{\sc Theorem  1.3.}
{\it
The image of the measure $\kappa_{n+1}$
under the map $\Upsilon$
coincides with $\kappa_n$.}

\smallskip

Hence, by the Kolmogorov theorem on the inverse limits,
there exists a canonical probabilitic measure
$\kappa_\infty$
on the space of infinite Rayleigh triangles
such that its projection to
 $\cR_n(-\infty,\infty)$
coincides with $\kappa_n$ for all $n$.

\smallskip

Possibly, this construction
is too general. We mention
its definitely reasonable partial case.
Assume
$\theta_{j\alpha}=\theta$
for all $j$, $\alpha$
and
$\sigma_j=\sigma+j\theta$.
Then we obtain a series of measures
interpolating the Hua-Pickrell measures
on the inverse limits of symmetric spaces
$$
\lim\limits_{\infty\leftarrow n}
\U(n)/\OO(n);\qquad
\lim\limits_{\infty\leftarrow n}
\U(n);\qquad
\lim\limits_{\infty\leftarrow n}
\U(2n)/\Sp(n)
;$$
these cases correspond to
$\theta=1/2,\,1,\, 2$.
With comparison with
0.3, we change the notation,
namely $\lambda_{j\alpha}$ is connected with $\nu_{j\alpha}$
from 1.3 by the formula
$\lambda_{j\alpha}=i (1-\nu_{j\alpha})/(1+\nu_{j\alpha})$.

\smallskip

{\sc Remark.} The space
$\cR_n[-\infty,\infty]$ discussed in this subsection
imitates the space of  Hermitian matrices,
 but in Subsection 0.3 we told about unitary matrices.
The passage from unitary matrices to Hermitian
matrices is given by the Cayley transform
$$
T=-1+(1+g)^{-1};\qquad g\in\U(n)
,\quad T\in\Herm_\C(n)
$$
For a matrix $T$,
we denote by $\langle T\rangle_j$
its lower right block of size $j\times j$.
Then (see \cite{Nerh})
$$
\det(1+[g]_{k} )=2^{-n+k}
\det(1+\langle T\rangle_{n-k}) \cdot \det(1+T)^{-1}
$$
and this allows to reduce the measures from
0.3 to the form
$c_n\cdot \fP_n(\cL) d\cL$.

\smallskip

{\bf 1.6. Trapezoids.}
Integrals (1.7)-(1.8), and also
integrals (3.1)--(3.3)
given below allow a minor extension.
For definiteness, we discuss
the integral (1.7).

Name by a  {\it Rayleigh trapezoid}
a collection of numbers
\def\aa{\scriptstyle}
$$           \aa
\begin{array} {ccccccccccccc}
 & & &\aa\lambda_{m1} &   &\aa\lambda_{m2}&\dots &  \aa\lambda_{mm} & & &  \\
 & &\aa\lambda_{(m+1)1} & &\aa \lambda_{(m+1)2}
           & &\dots  & &\aa\lambda_{(m+1)(m+1)} & &  \\
 &\cdot &\cdot & \cdot& \cdot &\cdot & \cdot &\cdot &\cdot &\cdot &  \\
\aa\lambda_{n1} & &\aa\lambda_{n2} & &\aa\lambda_{n3}
      & &  \dots & &\aa\lambda_{n(n-1)}
        &  & \aa\lambda_{nn}
,
\end{array}
,
$$
satisfying the interchanging condition (1.3).
We denote the set of all such trapezoids
by $\cR^m_n$.

The following formula holds
\begin{multline}
\int\limits_{\cR^m_n[0,\infty)}
\prod\limits_{j=m}^{n-1}\prod\limits_{\alpha=1}^j
\frac{\lambda_{j\alpha}^{\sigma_j-\sigma_{j+1}-\theta_{j\alpha}} }
 {(1+\lambda_{j\alpha})^{\tau_j-\tau_{j+1}+\theta_{j\alpha}} }
\cdot \prod\limits_{p=1}^n
  \frac {\lambda_{np}^{\sigma_n-1}}
     {(1+\lambda_{np})^{\tau_n}}
\times\\ \times
\prod\limits_{1\le \alpha<\beta\le m}
(\lambda_{m\beta}-\lambda_{m\alpha})^{2\kappa-1}
\times\\ \times
\frac{ \prod\limits_{j=m}^{n-1}
  \prod\limits_{1\le\alpha\le j,\, 1\le p\le j+1}
   |\lambda_{j\alpha}-\lambda_{(j+1)p}|^{\theta_{j\alpha}-1}
   }
   {   \prod\limits_{j=m}^{n-1}
  \prod\limits_{1\le\alpha<\beta \le j}
  (\lambda_{j\beta}-\lambda_{j\alpha})^{\theta_{j\alpha}+\theta_{j\beta}-2}
  }
\prod\limits_{1\le p< q\le n} (\lambda_{nq}-\lambda_{np})\, d\cL=
\\=
\prod\limits_{j=1}^{m-1}
\frac
{\Gamma((j+1)\kappa)\Gamma(\sigma_m+j\kappa)
        \Gamma(\tau_m-\sigma_m-(m+j-1)\kappa)}
{\Gamma(\kappa)\Gamma(\tau_m-j\kappa)}
\times\\ \times
\prod\limits_{m\le\alpha\le j\le n-1} \Gamma(\theta_{j\alpha})
\,\,\prod\limits_{j=m}^n
\frac
{\Gamma(\sigma_j)\Gamma(\tau_j-\sigma_j-
    \sum_{\alpha=1}^{j-1}\theta_{(j-1)\alpha})}
{\Gamma(\tau_j)}
.
\end{multline}

Emphasis, that for
 $m=n$ we obtain the Selberg integral
(see \cite{Sel}, \cite{And}, \cite{AAR}):
\begin{multline}
\int\limits_{0<\mu_1<\dots<\mu_n}
\prod\limits_{p=1}^n \frac{\mu_p^{\sigma-1}}{(1+\mu_p)^\tau}
\prod\limits_{1\le p<q\le n}
(\mu_q-\mu_p)^{2\theta}      d\mu_1\dots\mu_n =
\\
=\prod\limits_{j=1}^n
\frac{\Gamma(\sigma+\theta(j-1))
   \Gamma(\tau-\sigma-\theta(n-j+1))\Gamma(j\theta)}
   {\Gamma(\tau-\theta(j-1)\Gamma(\theta)}
\end{multline}
(we denote the variables
$\lambda_{np}$ by $\mu_p$ and $\kappa$
by $\theta$; the remaining parameters
for  $m=n$ are lacked).

\smallskip

{\bf 1.7. Conjecture.}
Our ultra-beta-integrals are some kind of a superstructure
over the Selberg integral. The Selberg integral
is the simplest representative of
multidimensional beta-integrals.
It seems likely that a superstructure of this type
exists for other beta-integrals
related to root systems of the series $A_n$,
for instance for the weight function of the Macdonald
polynomials. Anyway, for the Macdonald polynomials
are possible integral representations similar to
\cite{OO}, \cite{Kaz1}, see \cite{Kaz2}.

\smallskip

{\bf  1.8. Further structure of the paper.}
Our \S 2 contains main lemmas.
Using them,
in \S 3, we easily deduce the projectivity theorem,
the integrals (1.7)--(1.8), (1.9),
а также (3.1)--(3.3). In \S 4, Propositions 1.1-1.2 are proved,
in this section we discuss some more matrix integrals.

\medskip

{\large\bf \S 2. Change of variables}

\medskip

\addtocounter{sec}{1}
\setcounter{equation}{0}

We use changes of variables of the Anderson type
 \cite{And}.

\smallskip

{\bf 2.1. Basic lemmas.}
Fix real numbers
$a<b$; we allow
$a=-\infty$ or $b=\infty$.
Let fixed numbers $\lambda_1$, \dots $\lambda_{n-1}$
 satisfy the inequalities
\begin{equation}
a<\lambda_1<\dots< \lambda_{n-1}<b
.
\end{equation}
Denote by $\Xi[a,b;\lambda]$
the set of all collections $\mu_1$, \dots, $\mu_n$,
satisfying the inequality
\begin{equation}
 a<\mu_1<\lambda_1<\mu_2<\dots<\lambda_{n-1}<\mu_n<b
.
\end{equation}

{\sc Lemma 2.1.}
{\it Let complex numbers $\sigma$, $\tau$ и $\theta_1$,\dots, $\theta_{n-1}$
satisfy the conditions
$$
\Re \theta_\alpha>0, \qquad \Re \sigma>0;\qquad
\Re\tau>\Re\sigma -\sum_\alpha \Re \theta_\alpha
.$$
Then}
\begin{multline}
 \!\!\!\!\!\!
\int\limits_{\Xi[0,\infty;\lambda]}
\prod\limits_{p=1}^n \frac{\mu_p^{\sigma-1}}
                    {(1+\mu_p)^\tau}
\cdot
\prod\limits_{1\le\alpha\le n-1,\,1\le p\le n}
|\lambda_\alpha-\mu_p|^{\theta_\alpha-1}
\prod\limits_{1\le p <q\le n} (\mu_q-\mu_p)
\prod\limits_{1\le p\le n} d\mu_p
=\\=
\frac{\Gamma(\tau-\sigma-\sum \theta_\alpha)
\prod\limits_{1\le\alpha\le n-1}
\Gamma(\theta_\alpha)}
     {\Gamma(\tau)}
\times\\ \times
\prod\limits_{1\le\alpha<\beta\le n-1}
 (\lambda_\beta-\lambda_\alpha)^{\theta_\alpha+\theta_\beta-1}
\cdot
\prod\limits_{1\le \alpha\le n-1} \frac{\lambda_\alpha^{\sigma-1+\theta_\alpha}}
                       {(1+\lambda_\alpha)^{\tau-\theta_\alpha}}
.\end{multline}

{\sc Lemma 2.2.}
{\it Let complex numbers $\sigma$, $\tau$, $\theta_1$,\dots,$\theta_{n-1}$
satisfy the conditions
$$\Re \theta_\alpha>0;\quad
\Re ( \sigma +\tau-\sum_\alpha \Re \theta_\alpha -1)> 0
.$$
Then }
\begin{multline}
\int\limits_{\Xi[-\infty,\infty;\lambda]}
\prod\limits_{p=1}^n (1+i\mu_p)^{-\sigma}
                    (1-i\mu_p)^{-\tau}
\cdot
\prod\limits_{1\le\alpha\le n-1,\,1\le p\le n}
|\lambda_\alpha-\mu_p|^{\theta_\alpha-1}
\times \\ \times
\prod\limits_{1\le p <q\le n} (\mu_q-\mu_p)
\prod\limits_{1\le p\le n} d\mu_p
=\\=
\frac{\pi\, 2^{2-\sigma-\tau}\Gamma(\sigma+\tau-\sum \theta_\alpha-1)
\prod\limits_{1\le\alpha\le n-1}\Gamma(\theta_\alpha)}
     {\Gamma(\tau)\Gamma(\sigma)}
\times \\ \times
\prod\limits_{1\le\alpha<\beta\le n-1}
 (\lambda_\beta-\lambda_\alpha)^{\theta_\alpha+\theta_\beta-1}
\cdot
\prod\limits_{1\le \alpha} (1+i\lambda_\alpha)^{-\sigma-\theta_\alpha}
                       (1-i\lambda_\alpha)^{-\tau-\theta_\alpha}
.
\end{multline}

{\bf 2.2. Variants.}
There are three following variants of Lemma 2.1--2.2.

A)
Let $\Re \sigma>0$, $\Re \theta_\alpha>0$,
$\Re \psi>0$.
Then
\begin{multline}
\int\limits_{\Xi[0,\infty;\lambda]}
\prod\limits_{p=1}^n \mu_p^{\sigma-1}
e^{-\psi\mu_p}
\prod\limits_{1\le\alpha\le n-1,\,1\le p\le n}
|\lambda_\alpha-\mu_p|^{\theta_\alpha-1}
\times \\ \times
\prod\limits_{1\le p <q\le n} (\mu_q-\mu_p)
\prod\limits_{1\le p\le n} d\mu_p
=\\ =
\psi^{-\sigma-\sum \theta_\alpha}
\Gamma(\sigma)\prod_{\alpha=1}^{n-1}\Gamma(\theta_\alpha)
\prod\limits_{1\le\alpha\le n-1}
\lambda^{\sigma-1+\theta_\alpha}\,
e^{-\psi\lambda_\alpha}
\times \\ \times
\prod\limits_{1\le\alpha<\beta\le n-1}
 (\lambda_\beta-\lambda_\alpha)^{\theta_\alpha+\theta_\beta-1}
.\end{multline}

B)
Let $\Re \psi>0$, $\Re\theta_\alpha>0$.
Then
\begin{multline}
\int\limits_{\Xi[-\infty,\infty;\lambda]}
\prod\limits_{p=1}^n e^{-\psi \mu^2_p}
\prod\limits_{1\le\alpha\le n-1,\,1\le p\le n}
|\lambda_\alpha-\mu_p|^{\theta_\alpha-1}
\prod\limits_{1\le p <q\le n} (\mu_q-\mu_p)
\prod\limits_{1\le p\le n} d\mu_p
=\\
=         \sqrt{2\pi}
\exp\bigl\{-\frac\psi 2\sum_\alpha \lambda_\alpha^2\bigr\}
\psi^{-1/2-\sum \theta_\alpha}
\prod\limits_{\alpha=1}^{n-1}  \Gamma(\theta_\alpha)
\prod\limits_{1\le\alpha<\beta\le n-1}
 (\lambda_\beta-\lambda_\alpha)^{\theta_\alpha+\theta_\beta-1}
.
\end{multline}

C)
Let $-\infty<a<b<\infty$. Then
\begin{multline}
\int\limits_{\Xi[a,b;\lambda]}
\prod\limits_{p=1}^n
(\mu_p-a)^{\sigma-1} (b-\mu_p)^{\tau-1}
\times\\ \times
\prod\limits_{1\le\alpha\le n-1,\,1\le p\le n}
|\lambda_\alpha-\mu_p|^{\theta_\alpha-1}
\prod\limits_{1\le p <q\le n} (\mu_q-\mu_p)
\prod\limits_{1\le p\le n} d\mu_p
=\\ =
(b-a)^{\sigma+\tau+\sum \theta_\alpha-1}
\frac{\Gamma(\sigma)\Gamma(\tau)
   \prod\limits_{\alpha=1}^{n-1} \Gamma(\theta_\alpha)}
   {\Gamma(\sigma+\tau+\sum\theta_\alpha)}
 \times\\ \times
\prod_{\alpha=1}^{n-1}
(\lambda_\alpha-a)^{\sigma+\theta_\alpha-1}
      (b-\lambda_\alpha)^{\tau+\theta_\alpha-1}
\prod\limits_{1\le\alpha<\beta\le n-1}
 (\lambda_\beta-\lambda_\alpha)^{\theta_\alpha+\theta_\beta-1}
.
\end{multline}

\smallskip

{\bf 2.3.  Change of variables.}
Let the parameters $\lambda_1<\dots <\lambda_{n-1}$
be fixed,
and let the variables $\mu_1$, \dots, $\mu_n$
satisfy the interchanging conditions
\begin{equation}
 \mu_1<\lambda_1<\mu_2<\dots<\lambda_{n-1}<\mu_n
.
\end{equation}
Introduce the new variables
 $\xi_1$, \dots, $\xi_{n-1}$, $\eta$ by the formulas
\begin{gather}
\xi_\alpha=- \frac
   {\prod\limits_{1\le p\le n} (\mu_p-\lambda_\alpha)}
   {\prod\limits_{1\le \beta\le n-1, \beta\ne \alpha}
(\lambda_\beta-\lambda_\alpha)};\\
\eta=\sum \limits_{1\le p\le n} \mu_p-
 \sum\limits_{1\le \beta\le n-1}\lambda_\beta.
\end{gather}

{\sc Lemma 2.3.}
a) {\it The map
$(\mu_1, \dots, \mu_n)\mapsto (\xi_1,\dots,\xi_{n-1},\eta)$
given by formula {\rm (2.9)--(2.10)}
is a bijection of the domain
{\rm (2.8)}
to the domain }
\begin{equation}
\xi_1>0,\quad \xi_2>0,\quad \dots,\quad \xi_{n-1}>0,\quad \eta\in \R
.
\end{equation}

b) {\it If in addition $\lambda_1>0$,
then the inequalities
\begin{equation}
0<\mu_1< \lambda_1<\mu_2<\lambda_2<\dots
< \lambda_{n-1}<\mu_n
\end{equation}
 are equivalent to the inequalities}
\begin{equation}
\xi_1>0,\quad \xi_2>0,\quad\dots,\quad \xi_{n-1}>0,
\quad \eta
>\sum\limits_{\alpha=1}^{n-1}\frac{\xi_\alpha}{\lambda_\alpha}
.
\end{equation}

c) {\it The Jacobian of the change of variables
{\rm (2.9) -- (2.10)} is}
\begin{equation}
\frac{\prod\limits_{1\le p<  q \le n} (\mu_q-\mu_p)}
{\prod\limits_{1\le\alpha<\beta\le n-1}(\lambda_\beta-\lambda_\alpha)}
.
\end{equation}

{\sc Proof.}
Consider the function
\begin{equation}
R(x):=\frac{\prod\limits_{1\le p\le n} (x-\mu_p)}
        {\prod\limits_{1\le\beta\le n-1} (x-\lambda_\beta)}
.
\end{equation}
The residues of $R(x)$ are
\begin{equation}
\res\limits_{x=\lambda_\alpha} R(x)=-\xi_\alpha.
\end{equation}
Hence
\begin{equation}
R(x)=x-\eta-\sum\limits_{1\le\alpha\le n-1}
\frac {\xi_\alpha}{x-\lambda_\alpha}
.\end{equation}

Positivity  of the variables
 $\xi_\alpha$ immediately follows from the inequalities
(2.8).

Conversely, assume that all the
 $\xi_\alpha$ are positive.
The expansion (2.17) implies that
the function $R(x)$ is $+\infty$
at the left end of the interval
$(\lambda_\alpha, \lambda_{\alpha+1})$,
at the right end it equals $-\infty$. Hence
$R(x)$ has
at least one root on this interval.
The same is valid on the rays
$(-\infty,\lambda_1)$ and $(\lambda_{n-1},+\infty)$.
Hence, for positive $\xi_\alpha$
all the roots $\mu_p$ of the equation $R(x)=0$
are real and satisfy the inequalities (2.8).

Thus the statement a) is proved.
We know that the interval $(-\infty,\lambda_1)$
contains one root of the equation $R(x)=0$.
It is positive iff $R(0)<0$;
this proves b).

To evaluate  the Jacobian, we observe
$$
\frac{\partial \xi_\alpha}{\partial \mu_p}=
\frac  {\xi_\alpha} {\mu_p-\lambda_\alpha}
;\qquad
\frac{\partial\eta} {\partial \mu_p} =1
.
$$
Hence the Jacobian is
$$
\Delta\cdot\prod\limits_{\alpha=1}^{n-1}\xi_\alpha
,
$$
where
\begin{equation}
\Delta=
\det\begin{pmatrix}
\frac 1{\mu_1-\lambda_1} & \frac 1{\mu_1-\lambda_2} &\dots &
                  \frac 1{\mu_1-\lambda_{n-1}}& 1\\
\frac 1{\mu_2-\lambda_1} & \frac 1{\mu_2-\lambda_2} &\dots &
                  \frac 1{\mu_2-\lambda_{n-1}}& 1\\
 \vdots &\vdots &\ddots &\vdots &\vdots \\
\frac 1{\mu_n-\lambda_1} & \frac 1{\mu_n-\lambda_2} &\dots &
                  \frac 1{\mu_n-\lambda_{n-1}}& 1
\end{pmatrix}
.
\end{equation}

{\sc Lemma 2.4.}
{\it The following identity holds}
$$
\Delta=
\frac
{\prod\limits_{1\le p<  q \le n} (\mu_p-\mu_q)
\prod\limits_{1\le\alpha<\beta\le n-1}(\lambda_\beta-\lambda_\alpha)}
{\prod\limits_{1\le p\le n,\, 1\le\alpha\le n-1}(\mu_p-\lambda_\alpha)}
.$$

{\sc Proof.}
We use the Cauchy determinant
\begin{equation}
\det\limits_{1\le \alpha\le n,\, 1\le l\le n}
\frac 1 {\mu_\alpha-\lambda_l}=
\frac{\prod\limits_{1\le p<  q \le n} (\mu_q-\mu_p)
\prod\limits_{1\le\alpha<\beta\le n-1}(\lambda_\beta-\lambda_\alpha)}
{\prod\limits_{1\le p\le n,\, 1\le\alpha\le n}(\mu_p-\lambda_\alpha)}
.\end{equation}

 Further, we multiply the last column of the determinant by
 $\lambda_n$
and turn $\lambda_n$ to infinity.
This proves Lemma 2.4. \hfill $\square$

\smallskip

It remains to notice that
$$
|\xi_1 \xi_2\dots \xi_{n-1}|
=\frac
{\prod\limits_{1\le p\le n,\, 1\le\alpha\le n-1}|\mu_p-\lambda_\alpha|}
{\prod\limits_{1\le\alpha<\beta\le n-1}(\lambda_\beta-\lambda_\alpha)^2}
,$$
and this proves Lemma 2.3.
\hfill $\square$

\smallskip

{\sc Lemma  2.5.} {\it In notation of Lemma
{\rm 2.3}, for each $a\in \C$ the following identity
holds}
$$
\prod\limits_{1\le p\le n} (\mu_p+a)=
\prod\limits_{1\le \alpha\le n-1} (\lambda_\alpha+a)
\cdot \Bigl\{a+\eta-\sum\limits_{\alpha=1}^{n-1}
   \frac {\xi_\alpha} {\lambda_\alpha+a}\Bigr\}
.$$

{\sc Proof.}
Consider a function
$$Q(x)= \frac{\prod\limits_{1\le p\le n} (x-\mu_p)}
    {(x+a) \prod\limits_{1\le \beta\le n-1} (x-\lambda_\beta)}
.$$
Its residues are
$$
\res\limits_{x=\lambda_\alpha} Q(x)=
   \frac {-\xi_\alpha} {\lambda_\alpha+a}
;\qquad
\res\limits_{x=-a}     Q(x)       =
-\frac{\prod(\mu_p+a)}{\prod(\lambda_\beta+a)};
\qquad
\res\limits_{x=\infty}  Q(x) =\eta+a
.$$
The sum of residues is 0,
this is equivalent to the statement of Lemma.
        \hfill $\square$

\smallskip

{\bf 2.4. Deduction of Lemma 2.1.}
Denote our integral by $I_1$.
Let
\begin{equation}
W=
\prod\limits_{1\le\alpha<\beta\le n-1}
 (\lambda_\beta-\lambda_\alpha)^{\theta_\alpha+\theta_\beta-1}
.\end{equation}
Passing to the variables $\xi_\alpha$, $\eta$
using Lemma 2.3 и 2.5,
we obtain
\begin{multline}
\!\!\! \!\!\!\!\!\!
\frac {I_1} W=
\prod\limits_{\alpha=1}^{n-1}
\frac{\lambda_\alpha^{\sigma-1}}
  {(1+\lambda_\alpha)^{\tau}}
\int
\Bigl(\eta\,-\,\sum\frac{\xi_\alpha}{\lambda_\alpha}\Bigr)^{\sigma-1}
\Bigl(1+\eta\,-\sum\frac{\xi_\alpha}{1+\lambda_\alpha}\Bigr)^{-\tau}
\times \\ \times
\prod\limits_{\alpha=1}^{n-1} \xi_\alpha^{\theta_\alpha-1}
d\eta \prod\limits_{\alpha=1}^{n-1} d\xi_\alpha
,
\end{multline}
where the integration is given over the domain
$\xi_\alpha>0$,
$\eta>\sum\frac{\xi_\alpha}{\lambda_\alpha}$.

Further, we change $\eta$ и $\xi_\alpha$
to the new variables
\begin{gather*}
u=
\eta-\sum_{\alpha=1}^{n-1} \frac{\xi_\alpha}
                {\lambda_\alpha} ;\qquad\qquad
v_\alpha=\frac{\xi_\alpha}{\lambda_\alpha(1+\lambda_\alpha)}
.
\end{gather*}
In these variables,
$$
1+\eta-\sum_{\alpha=1}^{n-1} \frac{\xi_\alpha}
                {1+\lambda_\alpha}
=1+u+\sum_{\alpha=1}^{n-1} v_\alpha,
,$$
and
$$
\prod\limits_{\alpha=1}^{n-1}
\xi_\alpha^{-1}
\cdot d\eta\,d\xi_1\dots d\xi_{n-1}=
\prod\limits_{\alpha=1}^{n-1}
v_\alpha^{-1}
\cdot
du\,dv_1\dots dv_{n-1}
.
$$
The integral  converts to the form
\begin{multline*}
\prod\limits_{\alpha=1}^{n-1}
   \frac {\lambda_\alpha^{\sigma+\theta_\alpha-1}}
              {(1+\lambda_\alpha)^{\tau-\theta_\alpha}}
\int
\frac{u^{\sigma-1} \prod v_\alpha^{\theta_\alpha-1}}
{(1+u+v_1+\dots +v_{n-1})^\tau}
du\, dv_1\dots dv_{n-1}
,
\end{multline*}
where the integration is given over the domain
$u>0$, $v_\alpha>0$,
and we obtain a variant of the Dirichlet
integral, see the next subsection.

\smallskip

{\bf 2.5. Dirichlet integrals.}
The Dirichlet integrals are the three
following related integrals
\begin{gather}
\int\limits_{t_1>0,\dots, t_{n}>0,\sum t_n<1} \prod\limits_{j=1}^n
 t_j^{x_j-1} (1-\sum\limits_{j=1}^n t_j)^{x_{n+1}-1} dt_1\dots dt_{n}=
\frac{\prod\limits_{j=1}^{n+1} \Gamma(x_j)}{\Gamma(\sum\limits_{j=1}^{n+1}x_j) }
;\\
\int\limits_{t_1>0,\dots,t_n>0}
\frac{\prod\limits_{j=1}^n t_j^{x_j-1}dt_1\dots dt_n}
   {(1+\sum\limits_{j=1}^n  t_j)^y}=
\frac{\Gamma(y-\sum\limits_{j=1}^n x_j)\prod\limits_{j=1}^n \Gamma(x_j)}
       {\Gamma(y)}\qquad \qquad\qquad; \\
\int\limits_{t_1>0,\,\dots,\,t_{n-1}>0;\, s\in\R}
\frac{\prod\limits_{j=1}^{n-1} t_j^{x_j-1}ds\, dt_1\dots dt_{n-1}}
{(1+is+\sum\limits_{j=1}^{n-1} t_j)^y (1-is+\sum\limits_{j=1}^{n-1} t_j)^z}
=\qquad \qquad\qquad              \nonumber
;\\
\qquad \qquad\qquad\qquad\qquad\qquad
=\frac{2^{2-y-z} \pi \Gamma(y+z-1-\sum\limits_{j=1}^{n-1} x_j)
      \prod\limits_{j=1}^{n-1}\Gamma(x_j)}
  {\Gamma(y)\Gamma(z)}
.\end{gather}

All these  integrals can be easily evaluated
by successive integration in
$t_1$, $t_2$, \dots,
see also Theorem 1.8.6 from  \cite{AAR}.

\smallskip

{\bf 2.6. Deduction of Lemma 2.2.}
Denote our integral
by $I_2$, define the quantity $W$ by  (2.20).
Then by Lemma 2.3, 2.5 we obtain
\begin{multline*}
\frac {I_2} W=
\prod\limits_{\alpha=1}^{n-1} (1+i\lambda_\alpha)^{-\sigma}
 (1-i\lambda_\alpha)^{-\tau}
\times\\ \times
\int\Bigl(1+i\eta+\sum\limits_{\alpha=1}^{n-1}
  \frac{\xi_\alpha}{1+i\lambda_\alpha}\Bigr)^{-\sigma}
\Bigl(1-i\eta+\sum\limits_{\alpha=1}^{n-1}
  \frac{\xi_\alpha}{1-i\lambda_\alpha}\Bigr)^{-\tau}
\prod\limits_{\alpha=1}^{n-1} \xi_\alpha^{\theta_\alpha-1}
\,d\eta\,\prod d \xi_\alpha
,
\end{multline*}
where the integration is given over the domain
 (2.11).

The first and second factors of the integrand are converted to the form
$$
\Bigl(1+i\Bigl\{\eta-
 \sum\limits_{\alpha=1}^{n-1}
\frac{\xi_\alpha\lambda_\alpha}{1+\lambda_\alpha^2}
\Bigr\}
+\sum\limits_{\alpha=1}^{n-1}
  \frac{\xi_\alpha}{1+\lambda_\alpha^2}\Bigr)^{-\sigma}
\Bigl(1-i\Bigl\{
\eta-
 \sum\limits_{\alpha=1}^{n-1}
\frac{\xi_\alpha\lambda_\alpha}{1+\lambda_\alpha^2}
\Bigr\}
+\sum\limits_{\alpha=1}^{n-1}
  \frac{\xi_\alpha}{1+\lambda_\alpha^2}\Bigr)^{-\tau}
.$$
Further, introduce the new variables
$$
v_\alpha=\frac{\xi_\alpha}{1+\lambda_\alpha^2};
\qquad
u=
\eta-
 \sum\limits_{\alpha=1}^{n-1}
\frac{\xi_\alpha\lambda_\alpha}{1+\lambda_\alpha^2}
.
$$
Our integral transforms to the form
\begin{multline*}
\prod\limits_{\alpha=1}^{n-1} (1+i\lambda_\alpha)^{-\sigma+\theta_\alpha}
 (1-i\lambda_\alpha)^{-\tau+\theta_\alpha}
\times \\ \times
\int\limits_{v_1>0,\dots, v_{n-1}>0, u\in \R}
\frac{\prod v_\alpha^{\theta_\alpha-1} \, du\, \prod d v_\alpha\,\, }
{(1+iu+\sum v_\alpha)^\sigma  (1-iu+\sum v_\alpha)^\tau}
.
\end{multline*}
It remains to apply the Dirichlet integral
 (2.24).

\smallskip

{\bf 2.7. Deduction of integral (2.5).}
Denote our integral by $I_3$,
define the quantity $W$ by formula (2.20).
Then
$$
\frac {I_3} W=
\prod\limits_{\alpha=1}^{n-1}
\lambda_\alpha^{\sigma-1} e^{-\psi \lambda_\alpha}
\int \Bigl(\eta-\sum \frac{\xi_\alpha}{\lambda_\alpha}\Bigr)^{\sigma-1}
e^{-\psi\eta}\prod\limits_{\alpha=1}^{n-1} \xi_\alpha^{\theta_\alpha-1}
 \,d\eta\,  \prod\limits_{\alpha=1}^{n-1}  d \xi_\alpha
,
$$
where the integration is given over positive
 $\xi_\alpha$
and $\eta>\sum\xi_\alpha$.
Introducing the new variables
$$
u=\eta-\sum \frac{\xi_\alpha}{\lambda_\alpha}
;\qquad
v_\alpha=\frac{\xi_\alpha}{\lambda_\alpha}
,$$
we obtain
\begin{multline*}
\frac {I_3} W=
\prod\limits_{\alpha=1}^{n-1}
\lambda_\alpha^{\sigma+\theta_\alpha-1} e^{-\psi \lambda_\alpha}
\!\!\!
\int\limits_{u>0,v_1>0,\dots, v_{n-1}>0}
v^{\sigma-1}\prod v_\alpha^{\theta_\alpha-1}
\times\\ \times
\exp\bigl\{-\psi( u+\sum v_\alpha )\bigr\}
du\prod dv_\alpha
,
\end{multline*}
this comes to a product of
 $\Gamma$-functions.

\smallskip

{\bf 2.8. Derivation of formula (2.6).}

\smallskip

{\sc Lemma 2.6.}
{\it Let $\lambda_\alpha$, $\mu_p$, $\xi_\alpha$, $\eta$
be related by formula {\rm (2.9)--(2.10)}. Then}
$$\sum\limits_{p=1}^{n}\mu_p^2=
\sum\limits_{\alpha=1}^{n-1} \lambda_\alpha^2
+2\sum\limits_{\alpha=1}^{n-1} \xi_\alpha+\eta^2
.
$$

{\sc Proof.}
Represent the function
 (2.15) in the form
$$R(x)=x \prod \Bigl(1-\frac {\mu_p}x\Bigr)
         \prod \Bigl(1-\frac {\lambda_\alpha}x\Bigr)^{-1}
.
$$
The coefficient at
 $x^{-1}$ of its Laurent expansion at infinity
 is
$$
c=\sum\limits_{p<q} \mu_p\mu_q  +
 \sum\limits_{\alpha<\beta} \lambda_\alpha\lambda_\beta
-\sum\lambda_\alpha\mu_p+\sum\lambda_\alpha^2
=\frac12 \eta^2 +\frac 12 \sum\lambda_\alpha^2    -\frac12 \sum\mu_p^2
.
$$
On the other hand, the sum of residues is zero, i.e.,
 $c=-\sum\xi_\alpha$.
\hfill $\square$

\smallskip

Denote our integral (2.6) by $I_4$.
Passing to the variables $\xi_\alpha$, $\eta$,
we obtain the trivial integral
$$
\frac {I_4} W
=\exp\bigl\{-\frac\psi 2\sum\lambda_\alpha^2\bigr\}
\int_{\xi_\alpha>0, \eta\in\R}
\exp\bigl\{-\frac\psi 2 \eta^2-\psi\sum\xi_\alpha\bigr\}
\prod\xi_\alpha^{\theta_\alpha-1}d\eta\prod d\xi_\alpha
.$$

{\bf 2.9. Another deduction of integrals (2.5)--(2.6).}
Observe that
$$e^{-\mu}=\lim\limits_{k\to\infty}  \Bigl(1+\frac{\mu} k\Bigr)^{-k}
;\qquad
e^{-\mu^2}=\lim\limits_{k\to\infty}
  \Bigl(1+\frac{i\mu} k\Bigr)^{-k} \Bigl(1-\frac{i\mu} k\Bigr)^{-k}
,$$
hence the integral
 (2.5) can be obtained from (2.3)
by the substitution
$$\tau=k;\quad
\mu_p=\widetilde\mu_p/k; \lambda_\alpha=\widetilde\lambda_\alpha/k
$$
and passing to the limit as $k\to\infty$.
Similarly (2.6) can be derived from (2.4).

\smallskip

{\bf 2.10. Derivation of formula (2.7).}
In this case, it is possible to repeat literally
the calculation from Subsection 2.4 and
to reduce our integral to the Dirichlet integral
(2.23).

Also, it is possible to change the variables in (2.7)
in the following way
$$
\widetilde \mu_p =\frac{b-a}{\mu_p-a};\qquad
\widetilde \lambda_\alpha =\frac{b-a}{\lambda_\alpha-a}
$$
Thus we reduce our integral to
 (2.3).

\smallskip

{\bf 2.11. One more integral.}
Fix constants $\mu_1<\dots<\mu_{n+1}$.
Let the variables
$\lambda_1$, \dots $\lambda_{n}$
satisfy the inequalities
\begin{equation}
\mu_1<\lambda_1<\mu_2<\lambda_2<\dots<\lambda_{n}<\mu_{n+1}
.\end{equation}

{\sc Lemma 2.7.}
{\it Let $\Re\theta_1>0$, \dots, $\Re\theta_1>0$.
 Then}
\begin{multline}
\int \prod\limits_{1\le \alpha\le n;\, 1\le p\le n+1}
 |\lambda_\alpha-\mu_p|^{\theta_p-1}
\prod\limits_{1\le\alpha<\beta\le n}(\lambda_\beta-\lambda_\alpha)
 \prod \limits_{\alpha=1}^{n} d\lambda_\alpha=
\\ =
\prod\limits_{1\le p< q\le n+1}
(\mu_q-\mu_p)^{\theta_p+\theta_q-1}
\cdot
\frac{\prod\limits_{\alpha=1}^{n+1} \Gamma(\theta_\alpha)}
 {\Gamma(\sum\limits_{\alpha=1}^{n+1} \theta_\alpha)}
.
\end{multline}

{\bf 2.12. One more change of variables.}
Define the new variables $\zeta_1$,\dots, $\zeta_{n+1}$
by the formula

\begin{equation}
\zeta_p=\frac{\prod\limits_{1\le \alpha\le n} (\lambda_\alpha-\mu_p)}
   {\prod\limits_{1\le q\le n+1, \, q\ne p}(\,\, \mu_q-\mu_p)}
.\end{equation}

{\sc Lemma 2.8.} a)  {\it The map $\lambda\mapsto\zeta(\lambda)$
transforms the set defined by the inequalities {\rm (2.8)}
to the simplex}
$$.
\sum\limits_{p=1}^{n+1} \zeta_p=1; \qquad \zeta_1>0,\dots, \zeta_{n+1}>0
.$$

b) {\it  The Jacobian of the change of variables
$$(\lambda_1,\dots,\lambda_{n})\mapsto (\zeta_1, \dots, \zeta_{n})$$
is                  }
$$\frac
{\prod\limits_{1\le \alpha<\beta\le n}(\lambda_\beta-\lambda_\alpha)}
{\prod\limits_{1\le p<q\le n+1}(\mu_p-\mu_q)}
.
$$

{\sc Proof.}
Consider the function
$$
H(x)=\frac
    {\prod\limits_{1\le \alpha\le n} (x-\lambda_\alpha)}
{\prod\limits_{1\le p\le n+1} (x-\mu_p)}
.$$
The residues of this function are
$$\res\limits_{x=\mu_p}H(x)=\zeta_p;\qquad
 \res\limits_{x=\infty}H(x)=-1
.$$
The sum of all the residues is 0, hence $\sum\zeta_p=1$.
Obviously, inequalities (2.8) imply positivity of
$\zeta_p$.
Conversely, let all the $\zeta_p$ be positive.
Taking in account the expansion
$$
H(x)=\sum\limits_{p=1}^{n+1} \frac {\zeta_p}{x-\mu_p}
$$
we observe that the function
$H(x)$ has a unique zero on each interval
$(\mu_p,\mu_{p+1})$.
This proves statement a).

\smallskip

Further
$$
\frac {\partial \zeta_p}{\partial \lambda_\alpha}=
  \frac{\zeta_p}{\lambda_\alpha -\mu_p}
.
$$
Hence the evaluation of the Jacobian is reduced
to the Cauchy determinant
 (2.19).
\hfill $\square$

\smallskip

{\bf 2.13. Derivation of Lemma 2.7.}
The change of variables (2.27)
converts integral (2.26) to the form
$$
\prod\limits_{1\le p< q\le n+1}
(\mu_q-\mu_p)^{\theta_\alpha+\theta_\beta-1}
\int\limits_{\zeta_1>0,\dots,\zeta_n>0}
\prod\limits_{p=1}^{n+1}\zeta_p^{\theta_p-1}
\prod\limits_{p=1}^{n} d\zeta_p
.$$

Applying $\zeta_{n+1}=1-\zeta_1-\dots-\zeta_{n}$,
we obtain the Dirichlet integral (2.22).

\medskip

{\large\bf 3. Big integrals}

\medskip

\addtocounter{sec}{1}
\setcounter{equation}{0}

Here we prove Theorem 1.3 on projectivity,
and also deduce the ultra-beta-integrals
 (1.7)--(1.8),
(1.9), and similar integrals (3.1)--(3.3).

\smallskip

{\bf 3.1. Deduction of integrals (1.7)--(1.8).}
For definiteness, consider the integral (1.7).
Using Lemma 2.1 we integrate
in the variables
$\lambda_{n1}$, \dots,
$\lambda_{nn}$.
The result is an integral of the same form
but over the space $\cR_{n-1}[0,\infty)$
of Rayleigh triangles of lower size.

Next, we integrate in
$\lambda_{(n-1)1}$, \dots,
$\lambda_{(n-1)(n-1)}$, etc.

\smallskip

{\bf 3.2. Projectivity.}
In the notation of 1.1, let
$f(\cL)$ be a function on $\cR_{n+1}(-\infty,\infty)$,
depending only on the
first $n$ rows of a Rayleigh triangle.
To derive Theorem 1.3, we must prove that
$$
\int\limits_{\cR_{n+1}(-\infty,\infty)} f(\cL)\,d\kappa_{n+1}(\cL)=
\int\limits_{\cR_{n}(-\infty,\infty)} f(\cL)\,d\kappa_{n}(\cL)
$$

But Lemma 2.2 allows to integrate explicitly
over the variables
$\lambda_{(n+1)1}$, \dots,
$\lambda_{(n+1)(n+1)}$,
After this we obtain the required result.

\smallskip

{\bf 3.3. Variants.}
Formula (2.5)--(2.7) imply the following integrals
\begin{multline}
\int\limits_{\cR_n[0,\infty)}
\prod\limits_{j=1}^{n-1}\prod\limits_{\alpha=1}^j
\lambda_{j\alpha}^{\sigma_j-\sigma_{j+1}-\theta_{j\alpha}}
  e^{-(\psi_j-\psi_{j+1})\lambda_\alpha }
     \cdot \prod\limits_{p=1}^n
  \lambda_{np}^{\sigma_n-1}
     e^{-\psi_n \lambda_{np}}
\times\\ \times
\prod\limits_{j=1}^{n-1}
\frac{
  \prod\limits_{1\le\alpha\le j,\, 1\le p\le j+1}
   |\lambda_{j\alpha}-\lambda_{(j+1)p}|^{\theta_{j\alpha}-1}
   }
   {
  \prod\limits_{1\le\alpha<\beta \le j}
  (\lambda_{j\beta}-\lambda_{j\alpha})^{\theta_{j\alpha}+\theta_{j\beta}-2}
  }
\prod\limits_{1\le p< q\le n} (\lambda_{nq}-\lambda_{np})\, d \cL=
\\=
\prod\limits_{1\le\alpha\le j\le n-1} \Gamma(\theta_{j\alpha})
\prod\limits_{j=1}^n
  \Gamma(\sigma_j)
    \psi_j^{-\sigma_j-\sum\limits_{\alpha=1}^{n-1} \theta_{(j-1)\alpha}}
;\end{multline}

\begin{multline}
\int\limits_{\cR_n(-\infty,\infty)}
\prod\limits_{j=1}^{n-1}
\exp\Bigl\{-\frac12(\psi_j-\psi_{j+1})
   \sum\limits_{\alpha=1}^j\lambda_{j\alpha}^2\Bigr\} \cdot
\exp\Bigl\{-\frac12 \psi_n   \sum\limits_{p=1}^n \lambda_{np}^2\Bigr\}
\times\\ \times
\prod\limits_{j=1}^{n-1}
\frac{
  \prod\limits_{1\le\alpha\le j,\, 1\le p\le j+1}
   |\lambda_{j\alpha}-\lambda_{(j+1)p}|^{\theta_{j\alpha}-1}
   }
   {
  \prod\limits_{1\le\alpha<\beta \le j}
  (\lambda_{j\beta}-\lambda_{j\alpha})^{\theta_{j\alpha}+\theta_{j\beta}-2}
  }
\prod\limits_{1\le p< q\le n} (\lambda_{nq}-\lambda_{np})\, d \cL=
\\=
(2\pi)^{n/2}
\prod\limits_{1\le\alpha\le j\le n-1} \Gamma(\theta_{j\alpha})
\prod\limits_{j=1}^n
    \psi_j^{-n/2-\sum\limits_{\alpha=1}^{n-1} \theta_{(j-1)\alpha}}
;\end{multline}

\begin{multline}
\int\limits_{\cR_n[a,b]}
\prod_{j=1}^{n-1} \prod_{\alpha=1}^j
(\lambda_{j\alpha}-a)^{\sigma_j-\sigma_{j+1}-\theta_{j\alpha}}
(b-\lambda_{j\alpha})^{\tau_j-\tau_{j+1}-\theta_{j\alpha}}
\times\\ \times
\prod_{p=1}^n(\lambda_{np}-a)^{\sigma_n-1}
(b-\lambda_{np})^{\tau_n-1}
\times\\ \times
\prod\limits_{j=1}^{n-1}
\frac{
  \prod\limits_{1\le\alpha\le j,\, 1\le p\le j+1}
   |\lambda_{j\alpha}-\lambda_{(j+1)p}|^{\theta_{j\alpha}-1}
   }
   {
  \prod\limits_{1\le\alpha<\beta \le j}
  (\lambda_{j\beta}-\lambda_{j\alpha})^{\theta_{j\alpha}+\theta_{j\beta}-2}
  }
\prod\limits_{1\le p< q\le n} (\lambda_{nq}-\lambda_{np})\, d \cL=
\\=
(b-a)^{\sum\limits_{j=1}^n(\sigma_j+\tau_j-1)+
     \sum\limits_{1\le j\le\alpha\le n-1}\theta_{j\alpha}}
\times\\ \times
\prod\limits_{1\le j\le\alpha\le n-1} \Gamma(\theta_{j\alpha})
\cdot
\prod_{j=1}^n\frac{\Gamma(\sigma_j)\Gamma(\tau_j)}
                   {\Gamma(\sigma_j+\tau_j+
        \sum\limits_{\alpha=1}^{j-1} \theta_{(j-1)\alpha})}
.\end{multline}

{\bf 3.4. Trapezoids.}
We start from the integral
\begin{multline}
\int\limits_{\cR_n[0,\infty)}
\prod_{j=1}^{m-1}
\frac{\prod\limits_{1\le\alpha\le j;\,1\le p\le j+1}
   |\lambda_{j\alpha}-\lambda_{(j+1)p}|^{\kappa-1}}
   {\prod\limits_{1\le\alpha < \beta\le j}
  (\lambda_{j\beta}-\lambda_{j\beta})^{2\kappa-2}}
\times\\ \times
\prod_{j=m}^{n-1} \prod\limits_{\alpha=1}^j
 \frac{(\lambda_{j\alpha})^{\sigma_j-\sigma_{j+1}-\theta_{j\alpha}}}
 {(1+\lambda_{j\alpha})^{\tau_j-\tau_{j+1}+\theta_{j\alpha}}}
\cdot \prod_{p=1}^n \frac{\lambda_{np}^{\sigma_n-1}}
                   {(1+\lambda_{np})^{\tau_n}}
\times\\ \times
\prod_{j=m}^{n-1}
\frac{
  \prod\limits_{1\le\alpha\le j,\, 1\le p\le j+1}
   |\lambda_{j\alpha}-\lambda_{(j+1)p}|^{\theta_{j\alpha}-1}
   }
   {
  \prod\limits_{1\le\alpha<\beta \le j}
  (\lambda_{j\beta}-\lambda_{j\alpha})^{\theta_{j\alpha}+\theta_{j\beta}-2}
  }
\prod\limits_{1\le p< q\le n} (\lambda_{nq}-\lambda_{np})\, d \cL
,\end{multline}
which is a partial case of integral
(1.7) evaluated above.
The integration in $\lambda_{11}$ is reduced to the evaluation
of the integral
$$
\int\limits_{\lambda_{21}}^{\lambda_{22}}
(\lambda_{11}-\lambda_{21})^{\kappa-1}
(\lambda_{22}-\lambda_{11})^{\kappa-1}
    d\lambda_{11}
.
$$
After this integration, we have an  integral over
$\cR^2_n[0,\infty)$
of the form (1.9). Now the integration
in
$\lambda_{21}$, $\lambda_{22}$
is reduced to Lemma 2.7.
Further, using the same Lemma 2.7 we integrate the
last integral in
$\lambda_{31}$, $\lambda_{32}$, $\lambda_{33}$
etc.
As a result, we obtain the integral
(1.9) over  $\cR^m_n[0,\infty)$

\medskip

{\large\bf \S 4. Maps from matrices to Rayleigh triangles}

\medskip

\addtocounter{sec}{1}
\setcounter{equation}{0}

Here we prove Propositions 1.1--1.2,
and reduce some more matrix integrals
to the ultra-beta-integrals.

\smallskip

{\bf 4.1. Proof of Proposition 1.1.}
Consider the space $\Herm_\K(n)$
of Hermitian matrices $T$ of size $n\times n$
over $\K=\R,\C,\H$. Let us evaluate the distribution of
the eigenvalues $\mu_1$, \dots, $\mu_n$
of the matrix $T$ if
the block
 $[T]_{n-1}$
is fixed.

Write $T$ as a block $((n-1)+1) \times ((n-1)+1) $
matrix:
$$T=
\begin{pmatrix} T_{11} & T_{12} \\ T_{21} & T_{22}
\end{pmatrix}
.$$
By this time we fixed
 $T_{11}$; the block $T_{12}$
ranges in the coordinate space $\K^{n-1}$
with the standard Lebesgue measure;
the block $T_{22}$
ranges in real numbers,
at last, the block $T_{21}=T_{12}^*$
is completely defined by the block $T_{12}$.

Take a unitary matrix
 $g\in\U(n-1,\K)$,
such that $\Lambda=gT_{11}g^{-1}$
is a diagonal matrix.
Consider a matrix
$$
\begin{pmatrix} T'_{11} & T'_{12} \\ T'_{21} & T'_{22}
\end{pmatrix}
                                                 =
\begin{pmatrix} g&0\\0&1\end{pmatrix}
\begin{pmatrix} T_{11} & T_{12} \\ T_{21} & T_{22}
\end{pmatrix}
\begin{pmatrix} g^{-1}&0\\0&1\end{pmatrix}
=
\begin{pmatrix} \Lambda &g T_{12} \\ T_{21}g^{-1} & T_{22}
.
\end{pmatrix}
$$
We observe that the distribution of the vector $T'_{12}$
coincides with the distribution of the vector
$T_{12}$, а $T'_{22}=T_{22}$.

Thus, we reduce the problem  to the following form.
Let $T$ has the form
$$
T=\begin{pmatrix}
\lambda_1&0&\dots&0 & x_1\\
0&\lambda_2&\dots&0&  x_2\\
\vdots&\vdots&\ddots&\vdots&\vdots\\
0&0&\dots &\lambda_{n-1} &  x_{n-1}\\
\ov x_1& \ov x_2&\dots &\ov x_{n-1}&y
\end{pmatrix},
$$
where
$x_\alpha\in\K$, $y\in\R$. Let $\mu_1$,\dots, $\mu_n$
be its eigenvalues. Denote  by $\{dx\}$
the standard Lebesgue measure on the field $\K$ (see 1.2).
{\it We must find the image of the measure $\{dx_1\}\dots\{dx_{n-1}\}\,dy$
under the map }
\begin{equation}
(x_1,\dots,x_{n-1},y)\mapsto (\mu_1,\dots, \mu_n)
.
\end{equation}
Define the new variables
\begin{equation}
\xi_\alpha= |x_\alpha|;\qquad \eta=y
.\end{equation}
The image of the measure
$\{dx_1\}\dots\{dx_{n-1}\}dy$
under the map (4.2), obviously, is
\begin{equation}
\frac{\pi^{(n-1)\theta}}{\Gamma^{n-1}(\theta)}
\prod\limits_{\alpha=1}^{n-1} \xi_\alpha^{\theta-1} \,
\cdot \,d\eta  \prod\limits_{\alpha=1}^{n-1} d\xi_\alpha
.\end{equation}

The characteristic equation for the matrix
$T$ has the form
\begin{equation}
\mu-\eta-\sum_{\alpha=1}^{n-1} \frac{\xi_\alpha}{\mu-\lambda_\alpha}
=0
.
\end{equation}
Since the roots of this equation are
 $\mu_1$,\dots, $\mu_n$,
the left side of (4.4) can be rewritten in the form
$$
\frac
{\prod\limits_{1\le p\le n} (\mu-\mu_p)}
{\prod\limits_{1\le \beta\le n-1} (\mu-\lambda_\beta)}
.$$
From here, we obtain
\begin{gather*}
\xi_\alpha=- \frac
   {\prod\limits_{1\le p\le n} (\mu_p-\lambda_\alpha)}
   {\prod\limits_{1\le \beta\le n-1, \beta\ne \alpha}
(\lambda_\beta-\lambda_\alpha)},\\
\eta=\sum \limits_{1\le p\le n}\mu_p -
 \sum\limits_{1\le \beta\le n-1}\lambda_\beta
,
\end{gather*}
i.e., we obtain the familiar change of the variables
(2.9)--(2.10); in fact, this is the origin
of the formulas (2.9)--(2.10).

Applying Lemma 2.3  and formula (4.3),
we obtain that the distribution of the
eigenvalues $\mu_1$, \dots, $\mu_n$
is
$$
\frac{\pi^{(n-1)\theta}}{\Gamma^{n-1}(\theta)}
\cdot
\frac{\prod\limits_{1\le p\le n,\, 1\le\alpha\le n-1}
 |\mu_p-\lambda_\alpha|^{\theta-1}}
 {\prod\limits_{1\le\alpha<\beta\le n}
(\lambda_\beta-\lambda_\alpha)^{2\theta-2}}
  \cdot
\prod\limits_{1\le p<q\le n} (\mu_q-\mu_p)
\prod_{p=1}^n d\mu_p
.$$
This implies Proposition 1.1.

\smallskip

{\bf 4.2. Proof of Proposition 1.2.}
Consider the space
$\Mat_{n,n+k-1}(\K)$ of $n\times (n+k-1)$ matrices.

Represent such matrix as a block
$((n-1)+1)\times ((n-1)+k)$ matrix
$T=\begin{pmatrix} T_{11} & T_{12} \\ T_{21} & T_{22}
\end{pmatrix}$.
Evaluate distribution of the eigenvalues
of the matrix $TT^*$ if the
blocks $T_{11}$, $T_{12}$
are fixed.
For this, consider unitary matrices
$h\in \U(n-1,\K)$, $g\in \U(n-1+k)$,
such that the matrix
$h\begin{pmatrix} T_{11} & T_{12}
\end{pmatrix}g$
has the form $\begin{pmatrix} \Lambda & 0\end{pmatrix}$,
where $\Lambda$ is a diagonal matrix.
Consider the matrix
\begin{multline*}
\begin{pmatrix} T'_{11} & T'_{12} \\ T'_{21} & T'_{22}
\end{pmatrix}
:=
\begin{pmatrix} h&0\\0&1\end{pmatrix}
\begin{pmatrix} T_{11} & T_{12} \\ T_{21} & T_{22}
\end{pmatrix}
\begin{pmatrix} g_{11} & g_{12} \\ g_{21} & g_{22}
\end{pmatrix}
=\\=
\begin{pmatrix} \Lambda & 0\\
T_{21}\, g_{11}+T_{22} \, g_{21}&
T_{21}\, g_{12}+T_{22} \, g_{22}
\end{pmatrix}
\end{multline*}
We see that the block $T'_{12}$
is zero
and the vector-row $\begin{pmatrix} T_{21}'& T_{22}'\end{pmatrix}$
is distributed as
$\begin{pmatrix} T_{21}& T_{22}\end{pmatrix}$.

Thus the problem is reduced to the following  form.
Fix a matrix
$$
T=
\begin{pmatrix}
\lambda_1^{1/2}&0&\dots&0&0&\dots& 0\\
0&\lambda_2^{1/2}&\dots &0&0&\dots&0\\
\vdots &\vdots & \ddots &\vdots&\vdots& \ddots& \vdots&\\
0&0&\dots&\lambda_{n-1}^{1/2}&0&\dots&0\\
x_1&x_2& \dots &x_{n-1} & y_1& \dots & y_k
\end{pmatrix}
,$$
where $x_\alpha
$, $y_r\in\K$, and $\lambda_1$, \dots, $\lambda_{n-1}$
are real and nonnegative.
Let $\mu_1$, \dots, $\mu_n$ be
the eigenvalues of the matrix $TT^*$.
We must  {\it find the image of the Lebesgue measure
$\{dx_1\}\dots \{dx_{n-1}\} \{dy_1\}\dots \{dy_k\}$
under the map}
$$ (x_1,\dots, x_{n-1}; y_1,\dots, y_k)\mapsto
      (\mu_1,\dots,\mu_n)
.$$

Now
$$TT^*=
\begin{pmatrix}
\lambda_1&0&\dots&0& \lambda_1^{1/2} \overline x_1 \\
0&\lambda_2&\dots&0& \lambda_2^{1/2} \overline  x_2\\
\vdots& \vdots& \ddots&\vdots&\vdots \\
0&0&\dots&\lambda_{n-1}&\lambda_{n-1}^{1/2} \overline  x_{n-1} \\
\lambda_1^{1/2} x_1&\lambda_2^{1/2} x_2& \dots &\lambda_{n-1}^{1/2}x_{n-1} &
  \sum |x_\alpha|^2+\sum |y_l|^2
.\end{pmatrix}
$$
Define the new variables
$$
\xi_\alpha=|x_\alpha|^2;\qquad \eta=\sum\limits_{l=1}^k |y_l|^2
.$$
    The image of the measure
$\{dx_1\}\dots \{dx_{n-1}\} \{dy_1\}\dots \{dx_k\}$
          under the map
$$ (x_1,\dots, x_{n-1}; y_1,\dots, y_k)\mapsto
  (\xi_1,\dots,\xi_{n-1}, \eta)
,
$$
obviously is
$$
\frac{\pi^{(n-1)\theta}} {\Gamma^{n-1}(\theta)}
\frac{\pi^{k\theta}}{\Gamma(k\theta)}
\eta^{k\theta-1}
\cdot\prod\limits_{k=1}^n \xi_\alpha^{\theta-1}
 \,d\eta\,\prod d\xi_\alpha
.$$

Now, the equation on the eigenvalues
 $\mu$ of the matrix $TT^*$ has the form
\begin{equation}
\eta+\sum\limits_{\alpha=1}^{n-1}\xi_\alpha-\mu -\sum_{\alpha=1}^{n-1}
\frac{\lambda_\alpha \xi_\alpha}
{\lambda_\alpha-\mu}=0
.\end{equation}

Let $\mu_1$, \dots, $\mu_n$ ---
be the roots of this equation;
in other words we represent the left part
of  equation (4.6) in the form
$\frac{\prod(\mu_p-\mu)}{\prod(\lambda_\alpha-\mu)}$. Obviously,
\begin{align}
\xi_\alpha=- \frac 1{\lambda_\alpha}  \frac
   {\prod\limits_{1\le p\le n} (\mu_p-\lambda_\alpha)}
   {\prod\limits_{1\le \beta\le n-1, \beta\ne \alpha}
(\lambda_\beta-\lambda_\alpha)};\\
\eta+\sum\xi_\alpha=\sum \mu_p - \sum \lambda_\alpha
.
\end{align}

{\sc Lemma 4.1.}
{\it  The following identity holds}
\begin{equation}
\eta=\frac{\prod\limits_{p=1}^n\mu_p}
  {\prod\limits_{\beta=1}^{n-1} \lambda_\beta}
.\end{equation}

{\sc Proof.}
Consider the function
$$
U(x)=\frac{\prod(x-\mu_p)}
 {x\prod(x-\lambda_\alpha)}
.$$
Its residues are
$$\res\limits_{x=\lambda_\alpha}U(x)=-\xi_\alpha;\qquad
\res\limits_{x=0}U(x)=
-\frac{\prod\mu_p}{\prod \lambda_\alpha};\qquad
\res\limits_{x=\infty}U(x)=\sum \mu_p-\sum\lambda_\alpha
.
$$
The sum of residues is 0. Comparing this with (4.8),
we obtain (4.9). \hfill $\square$

\smallskip

Thus, $\xi_\alpha$ and $\eta$ are expressed in terms of
$\mu_1,\dots, \mu_p$ by formulas (4.7), (4.9).
The Jacobian of the transformation
$(\mu_1,\dots, \mu_n)\mapsto (\xi_1,\dots,\xi_{n-1},\eta)$
can be easily evaluated
 using the Cauchy determinant (2.19).
As a result, we obtain
\begin{equation}
d\xi_1\dots d\xi_{n-1}d\eta=
\frac{\prod\limits_{1\le p<q\le n} (\mu_q-\mu_p)}
  {\prod\limits_{1\le\alpha \le n-1}  \lambda_\alpha
\prod\limits_{1\le\alpha <\beta\le n-1}
(\lambda_\beta-\lambda_\alpha)} d\mu_1\dots d\mu_n
.\end{equation}

It remains to substitute
(4.7), (4.9) and (4.10) to   (4.5).

\smallskip

{\bf 4.3. Matrix integrals.} Let  $m\le n$.
Denote by ${\rm B}_{n,m}(\K)$
the set of all matrices over $\K=\R$, $\C$, $\H$
of  size
$n\times m$ with norm  $<1$;
the term 'norm' means the norm of an operator
from the Euclidean space
$\K^n$ to the Euclidean space $\K^m$.
As above, let  $\theta=\dim \K/2$.

Integral (3.3) and Proposition 1.2
obviously imply

\begin{multline}
\int\limits_{{\rm B}_{n,m}(\K)}
\prod\limits_{j=1}^{n-1}
\Bigl[
\det\bigl(\{T\}_j \{T\}_j^*\bigr)^{\sigma_j-\sigma_{j+1}-\theta}
\det\bigl(1-\{T\}_j \{T\}_j^*\bigr)^{\tau_j-\tau_{j+1}-\theta}
\Bigr]
\times\\ \times
\det(TT^*)^{\sigma_n-(m-n+1)\theta} \det(1-TT^*)^{\tau_n-1} dT
=\\=
\pi^{mn\theta}
\prod_{j=1}^n \frac{\Gamma(\sigma_j)\Gamma(\tau_j)}
            {\Gamma(\sigma_j+\tau_j+(j-1)\theta)\Gamma((m-j+1)\theta)}
.\end{multline}

Another partial case  of integral
 (3.3) is
(see Proposition 1.1) is the Gindikin beta-function
 (see \cite{Sie}, \cite{Gin},\cite{FK})
\begin{multline}
\int\limits_{T=T^*,\, 0<T<1}
 \prod_{j=1}^{n-1} \det[T]_j^{\sigma_j-\sigma_{j+1}-\theta}
\det(1-[T]_j)^{\tau_j-\tau_{j+1}-\theta}
\times\\ \times
\det T^{\sigma_n-1} \det(1-T)^{\tau_n-1}dT       =
\pi^{n(n-1)\theta}
\prod_{j=1}^n \frac{\Gamma(\sigma_j)\Gamma(\tau_j)}
            {\Gamma(\sigma_j+\tau_j+(j-1)\theta)}
,\end{multline}
where the integration is given over the set
of all the Hermitian matrices
over $\K$ satisfying the condition $0<T<1$.

The Whishart--Gindikin matrix gamma-function
\begin{multline}
\int\limits_{\Pos_\K(n)}
\prod\limits_{j=1}^{n-1}\det[T]_j^{\sigma_j-\sigma_{j+1}-\theta}
\exp\bigl\{-(\psi_j-\psi_{j+1})\tr[T]_j\bigr\}
\times\\ \times
\det T^{\sigma_n-1}
\exp\bigl\{-\psi_n\tr T\bigr\}dt
=
\pi^{n(n-1)\theta}
\prod_{j=1}^n \Gamma(\sigma_j) \psi_j^{-\sigma_j-(j-1)\theta}
\end{multline}
is a partial case of the integral
 (3.1).
Another partial case of
 (3.1) is the integral
\begin{multline}
\int\limits_{\Mat_\K(n,m)}
\prod\limits_{j=1}^{n-1}
\det(\{T\}_j\{T\}_j^*)^{\sigma_j-\sigma_{j+1}-\theta}
\exp\bigl\{-(\psi_j-\psi_{j+1})\tr  \{T\}_j\{T\}_j^* \bigr\}
\times \\ \times
\det(TT^*)^{\sigma_n-\theta(n-m+1)}
\exp\bigl\{-\psi_n\tr TT^*\bigr\} dT
=\\=
\pi^{mn\theta}
\prod_{j=1}^n
\frac{\Gamma(\sigma_j) \psi_j^{-\sigma_j-(j-1)\theta}}
{\Gamma((m-j+1)\theta)}
\end{multline}

The Mehta type integral (3.3) interpolates the  matrix
integrals
$$
\int\limits_{\Herm_{\K}(n)}
\exp\Bigl\{-\frac12\sum\limits_{j=1}^{n-1} (\psi_j-\psi_{j+1})\tr [T]_j^2
-\frac12\psi_n \tr T^2\Bigr\} \,dT
.
$$
However these integrals themselves are quite obvious,
since the expression in the exponent is a diagonal
quadratic form.

\sf
Institute of the Theoretical and Experimental Physics,

B. Cheremushkinskaya, 25

Moscow 117259

Russia


{\tt neretin@mccme.ru}

\end{document}